\documentclass[oneside,notitlepage,12pt]{article}

\usepackage{amssymb}
\usepackage[leqno]{amsmath}
\usepackage{amsfonts}
\usepackage{amsopn}
\usepackage{amstext}
\usepackage{amsthm}

\usepackage[all]{xy}

\usepackage[colorlinks,backref]{hyperref}
\usepackage{makeidx}

\usepackage{calrsfs}
\usepackage{fourier-orns}
\usepackage{enumitem}
\usepackage{clock} 
\ClockFrametrue\ClockStyle1

\newenvironment{pf}{\begin{proof}}{\end{proof}}

\newcommand{\El}{{\cal{L}}}

\newcommand{\Nat}{{\mathbb{N}}}

\newcommand{\sig}{\sigma}
\newcommand{\eps}{\varepsilon}
\renewcommand{\phi}{\varphi}
\renewcommand{\rho}{\varrho}

\newcommand{\ntr}{n\in\omega}

\newcommand{\loe}{\leqslant}
\newcommand{\goe}{\geqslant}

\newcommand{\subs}{\subseteq}

\newcommand{\nnempty}{\ne\emptyset}

\newcommand{\ovr}{\overline}
\renewcommand{\iff}{\Longleftrightarrow}

\newcommand{\id}[1]{{\operatorname{i\!d}_{#1}}} 

\newcommand{\dom}{\operatorname{dom}}

\newcommand{\cod}{\operatorname{cod}}
\newcommand{\Aut}{\operatorname{Aut}}

\newcommand{\oraz}{\qquad\text{and}\qquad}
\newcommand{\mion}[1]{{#1}^{-1}}
\newcommand{\Lip}[1]{\operatorname{Lip}\left(#1\right)}

\newcommand{\ob}[1]{\operatorname{Obj}\left(#1\right)}

\newcommand{\metric}{\mathfrak M\!\mathfrak s}

\newcommand{\epm}{\metric^+_{\rm{ep}}}
\newcommand{\nls}{{\mathfrak V}_{\anorm}}

\newtheorem{tw}{Theorem}[section]
\newtheorem{wn}[tw]{Corollary}
\newtheorem{lm}[tw]{Lemma}
\newtheorem{prop}[tw]{Proposition}
\newtheorem{claim}[tw]{Claim}

\theoremstyle{definition}
\newtheorem{df}[tw]{Definition}

\newtheorem{ex}[tw]{Example}

\theoremstyle{remark}

\newcommand{\Fr}{\operatorname{Fr}}

\providecommand{\nat}{\omega}

\newcommand{\setof}[2]{\{#1\colon #2\}}
\newcommand{\bigsetof}[2]{\Bigl\{#1\colon #2\Bigr\}}
\newcommand{\seq}[1]{\langle #1 \rangle}

\newcommand{\sett}[2]{\{#1\}_{#2}}
\newcommand{\sn}[1]{\{#1\}} 
\newcommand{\dn}[2]{\{#1,#2\}} 
\newcommand{\pair}[2]{{\left\langle{#1},{#2}\right\rangle}} 
\newcommand{\triple}[3]{\langle #1, #2, #3 \rangle} 
\newcommand{\map}[3]{#1\colon #2 \to #3} 
\newcommand{\inv}[2]{{#1}^{-1}[#2]} 

\newcommand{\ciag}[1]{{\sett{{#1}_n}{\ntr}}}

\newcommand{\iso}{\approx}

\newcommand{\anorm}{\|\cdot\|}
\newcommand{\norm}[1]{\|#1\|}

\newcommand{\fK}{{\mathfrak{K}}}
\newcommand{\fL}{{\mathfrak{L}}}

\newcommand{\fra}{Fra\"iss\'e}

\newcommand{\cmp}{\circ} 

\newcommand{\wek}[1]{{\vec{#1}}}

\newcommand{\separator}{\begin{center} \leafright \leafright \decotwo \leafleft \leafleft \end{center}}

\newcommand{\define}[2]{{\em #1}\index{#2}}

\newcommand{\ciagi}[1]{\sig{#1}}

\newcommand{\MCsymb}{\mathcal{N}}
\newcommand{\MC}[1]{{\ensuremath{\MCsymb_{#1}}}}
\newcommand{\MCf}{{\ensuremath{\MCsymb^+_3}}}

\title{Categories with norms}
\author{
Wies{\l}aw Kubi\'s
\footnote{Research supported by the GA\v{C}R project No. 17-27844S and RVO: 67985840.}
\\ {\small Institute of Mathematics}\\
{\small Czech Academy of Sciences, Czechia}\\
-----\\
{\small Institute of Mathematics}\\
{\small Cardinal Stefan Wyszy\'nski University in Warsaw, Poland}
}

\date{\clocktime\ \today}

\makeindex

\begin{document}

\maketitle

\begin{abstract}
We study a metric-like structure on categories, showing that the concept of the limit of a sequence in a metric space and the concept of the colimit of a sequence in a category have a common generalization. The main concept is a \emph{norm} on a category, generalizing pseudo-metrics and group valuations.
In this new context, we discuss topics like Cauchy completion and the Banach Contraction Principle.


\end{abstract}

\tableofcontents

\section{Introduction}

A celebrated work of Lawvere~\cite{Lawvere} exhibits surprising connections between metric spaces and categories.
A crucial link is the idea of enriched category, where several different concepts can be uniformized.

The aim of our note is to generalize the concept of a \emph{metric} in a different way, namely, where the underlying set is replaced by a category.
Our starting point is the observation that every set can be regarded as a (quasi-ordered) category, in which for every two objects $a$, $b$ there exists exactly one arrow $ab$ from $a$ to $b$.
The distance between $a$ and $b$ can be used as some ``measure" of the arrow $ab$. Since $ab$ also resembles a vector, its measure could actually be called 
\emph{length} or \emph{norm}.
This leads to the concept of a \emph{normed category}, where each arrow $f$ has a non-negative value $\mu(f)$ (possibly $+\infty$) and three natural axioms are imposed.
We give the precise definition in Section~\ref{SecMetsCatsERTrf}.
Presenting several examples, we argue that our list of axioms for a norm on a category is natural.
In case where the category is a quasi-ordered set in which each two elements are related, our axioms reduce to the usual axioms of a metric (including symmetry!) except that we allow the distance to be $+\infty$ and we also allow that different points might have distance zero (thus, formally, we obtain a metric after a suitable identification of points).
Our motivating example of a normed category is the category of metric spaces with Lipschitz mappings, where $\mu(f)$ is the logarithm of the maximum of the Lipschitz constants $\Lip{f}$, $\Lip{\mion f}$, agreeing that $\mu(f) = +\infty$ if either $f$ is not one-to-one or $\mion f$ is not Lipschitz.
This example also shows that allowing $\mu(f) = +\infty$ indeed makes sense, as this carries information that $f$ is not bi-Lipschitz.

Section~\ref{SecMetsCatsERTrf} contains, besides the crucial definition of a norm on a category, several natural examples of normed categories.
We show that even discrete norms (having only one positive value) lead to non-trivial examples, contrary to metric spaces, where a discrete metric is totally determined by the cardinality of the set of points.
We also show that asymmetric spaces can be viewed as normed categories, however one has to use certain directed graphs instead of quasi-orderings.

Section~\ref{SectCauchySeqsCmps} introduces the concept of a \emph{Cauchy sequence} and of its limit, which is the colimit in category-theoretic sense, requiring that the colimiting cocone of arrows satisfies a simple and natural convergence condition.
This is the place where the classical concept of a convergent sequence in a metric space meets the concept of a convergent sequence in a category.
Recall that a \emph{sequence} in a category is simply a covariant functor from the set of nonnegative integers $\nat$ treated as a poset category.
We shall describe a rather simple construction of the \emph{Cauchy completion}, namely, a bigger normed category in which every Cauchy sequence is convergent.

Section~\ref{SectBanachFixedpoint} contains a short discussion of Banach's Fixed-point Principle in normed categories.
It turns out that the main ideas of Banach's Principle can be easily rephrased in our setting, although the existence of a fixed-point (or rather a fixed object) is subject to some extra condition, which is always satisfied in a metric space.
Uniqueness of fixed-points is not a general phenomenon in our setting.


\section{Normed categories} \label{SecMetsCatsERTrf}

We shall use standard notation concerning category theory, see e.g. \cite{MacLane}.
Given a category $\fK$, we shall denote by $\ob \fK$ its class of objects and $\fK$ itself will serve as the class of arrows.
In particular, given $f \in \fK$, we shall denote by $\dom(f)$ and $\cod(f)$ its domain and codomain.
Given objects $x,y \in \ob \fK$, we shall denote by $\fK(x,y)$ the set of all $f \in \fK$ with $\dom(f) = x$, $\cod(f) = y$.
Obviously, $\map f x y$ is a shortcut for $f \in \fK(x,y)$.
We shall sometimes write ``$f$ is a $\fK$-arrow" instead of ``$f\in \fK$".
The composition operation in a category will always be denoted by $\cmp$.
More precisely, if $f \in \fK(y,z)$ and $g \in \fK(x,y)$ then $f \cmp g \in \fK(x,z)$.
The identity of $x \in \ob \fK$ will be denoted by $\id x$.
Below is the crucial definition.

\begin{df}
	A \define{norm}{norm} on a category $\fK$ is a function $\map \mu \fK{[0,+\infty]}$ satisfying the following conditions:
	\begin{enumerate}[itemsep=0pt]
		\item[(\MC1)] $\mu(\id x)=0$ for every object $x$.
		\item[(\MC2)] $\mu(f\cmp g) \loe \mu(f) + \mu(g)$ whenever $f\cmp g$ is defined.
		\item[(\MC3)] $\mu(g) \loe \mu(f\cmp g) + \mu(f)$ whenever $f\cmp g$ is defined.
	\end{enumerate}
	A pair $\pair \fK \mu$ will be called a \define{normed category}{category!-- normed}.
	\index{normed category}
\end{df}
Condition (\MC2) may be called, as usual, \define{triangle inequality}{triangle inequality}.
It seems reasonable to allow $\mu(f)=+\infty$, although it is always possible to replace $\mu(f)$ by $\mu_1(f) = \min(1,\mu(f))$, obtaining a new norm that is bounded by $1$.
The number $\mu(f)$ is meant to norm the `distortion' from the domain of $f$ to its range, or rather how $f$ transforms its domain into its codomain.
Perhaps condition (\MC3) requires some explanations.
First of all, note one of its consequences:

\begin{lm}\label{Cwergow}
$\mu(h) = \mu(h^{-1})$ whenever $h$ is an isomorphism.
\end{lm}

\begin{pf}
Let $x = \dom(h)$.
By (\MC3) and (\MC1), $\mu(h) \loe \mu(\id x) + \mu(h^{-1}) = \mu(h^{-1})$.
Replacing $h$ by $h^{-1}$ leads to $\mu(h^{-1}) \loe \mu(h)$.
\end{pf}

We say that $\map h a b$ is a \define{$0$-isomorphism}{$0$-isomorphism} if $h$ is an isomorphism and $\mu(h) = 0$.
In that case, also $\mu(h^{-1}) = 0$.
It is natural to define
$$\fK_0 = \setof{f \in \fK}{\mu(f) = 0}.$$
We shall call $\fK_0$ the \define{kernel}{normed category!-- kernel} of $\pair \fK \mu$.
Conditions (\MC1), (\MC2) say that $\fK_0$ is indeed a subcategory of $\fK$, while Lemma~\ref{Cwergow} implies that $\fK$-isomorphisms which are members of $\fK_0$ remain isomorphisms in $\fK_0$.
Actually, $\fK_0$ satisfies:
\begin{enumerate}
\item[(K1)] $\ob {\fK_0} = \ob \fK$,
\item[(K2)] $f \in \fK_0$ and $f \cmp g \in \fK_0 \implies g \in \fK_0$.
\end{enumerate}
We shall see later (Subsection~\ref{SubSectDiscreteCase} below) that every pair of categories $\fK_0 \subs \fK$ satisfying (K1), (K2) gives rise to a discrete (taking two values only) norm on $\fK$ such that $\fK_0$ becomes its kernel.

Coming back to condition (\MC3), let us note that if $\mu(f)$ norms some `distortion' of the codomain of $f$ with respect to its domain, then (\MC3) says that $g$ cannot make more distortion than the sum of the distortions of $f\cmp g$ and $f$.
As an example, let $\fK$ be the category of metric spaces and the distortion of a mapping $f$ is normed by the Lipschitz constants of $f$ and its inverse $f^{-1}$ (that is why we need the value $+\infty$ in case $f$ is not one-to-one).
The details are explained in Subsection~\ref{SubSectMetsNMesur}.
It seems that omitting (\MC3) would cause more `negative' consequences than omitting symmetry in the definition of a metric, see Subsection~\ref{SubSectMetricSpacesInGeneral}.
On the other hand, one may wonder why we do not require a stronger version of (\MC3), namely:
\begin{enumerate}
\item[(\MCf)] $|\mu(f)-\mu(g)| \loe \mu(f \cmp g)$ whenever $f\cmp g$ is defined.
\end{enumerate}
It turns out that this would be too strong, an evidence is given in the list of several natural examples of norms collected below, even though few of them actually satisfy (\MCf).

\subsection{Historical comments}

We feel that some explanations concerning our terminology are needed, supported by the (known to us) literature.

First of all, the name \emph{norm} typically refers to vector spaces, meaning the \emph{length} of their vectors. A vector actually represents a unique way of `moving' from its origin to its destination (end-point).
A vector space, as any semigroup, is indeed a category in which the composition is simply the addition of vectors. Scalar multiplications are ignored here. Anyway, it is clear that every norm on a real vector space is a norm in the sense of our definition, whenever the vector space is treated as a category.

Having in mind normed linear spaces and bounded linear operators, one immediately finds out that it is natural to consider categories enriched over normed spaces, namely, assuming that the hom-sets are normed vector spaces, requiring that the composition is bi-linear, and imposing some further conditions. By this way, the norm of an arrow is just the norm in the usual sense. This, unfortunately, formally differs from our setting, as $\norm{f \cmp g} \loe \norm{f} \cdot \norm{g}$ and a logarithm is required to get (\MC2) and (\MC3), see Subsection~\ref{SubSectBanachNorms} for details. Fortunately, this is just a formal (arithmetic) difference. We shall never use $\anorm$ for denoting a category norm.

The name \emph{norm} in the context of category theory, as one could expect, has already appeared in the literature (see \cite{Kent, Zlatos, Macaj, Grandis}), usually related to theoretical computer science.
In all instances, condition (\MC3) is discarded. In \cite{Zlatos, Macaj} the name \emph{valuation} was used instead of a norm, declaring clauses (\MC1), (\MC2), not (\MC3).
On the other hand, a valuation satisfying (\MC1)--(\MC3) appeared and played a significant role already in~\cite{America}, where the authors aimed at certain fixed-point theorems solving some domain equations.
We present this valuation in Subsection~\ref{SubSectionAMericaForEver}, showing that it is indeed a norm.
Last but not least, our motivation for studying norms on categories comes from the theory of \fra\ limits in a continuous setting, see Section~\ref{SectFinals} and our preprint~\cite{Kub41}.

Summarizing, we have decided to use the name `norm' for a function on arrows satisfying (\MC1)--(\MC3), taking into account the fact that such a name has already appeared in the category-theoretic literature. We hope that our examples support the claim that the axioms of a norm we propose are quite natural and applicable, especially when adding the somewhat controversial condition (\MC3).

\subsection{Discrete norms}\label{SubSectDiscreteCase}

Let us call a norm $\mu$ on a category $\fK$ \define{discrete}{norm!-- discrete} if it takes values in the set $\dn 0 {+\infty}$.
This is motivated by the notion of a discrete metric, where typically the values are in the set $\dn 01$ instead of $\dn 0{+\infty}$, which does not make any difference on the structure of a space.
While discrete metric spaces are not particularly interesting (precisely: the cardinality determines their isometric type), we shall see that discrete norms lead to interesting pairs of categories.

\begin{prop}
Let $\fK_0 \subs \fK$ be a pair of categories.
The following properties are equivalent:
\begin{enumerate}
\item[{\rm(a)}] There exists a norm on $\fK$ such that $\fK_0$ is its kernel.
\item[{\rm(b)}] There exists a discrete norm on $\fK$ with $\fK_0$ being its kernel.
\item[{\rm(c)}] $\fK_0$ satisfies
\begin{itemize}
\item[{\rm(K1)}] $\ob{\fK_0}=\ob\fK$,
\item[{\rm(K2)}] $\dn f {f \cmp g} \subs \fK_0 \implies g \in \fK_0$.
\end{itemize}
\end{enumerate}
\end{prop}

\begin{pf}
It is clear that (b)$\implies$(a) and we have already seen that (a)$\implies$(c) (which immediately follows from (\MC1)--(\MC3)).
It remains to show that (c)$\implies$(b).
Namely, we need to check that the function $\map \mu \fK {\dn 0 {+\infty}}$, defined by $\mu(f)=0$ if and only if $f \in \fK_0$, is a norm on $\fK$.

Condition (\MC1) follows from the fact that $\fK_0$ has all the objects of $\fK$, which is declared in (K1).
Condition (\MC2) follows from the fact that $\fK_0$ is a subcategory of $\fK$.
Finally, (\MC3) follows directly from (K2).
\end{pf}

Let us say that a subcategory $\fK_0$ of a category $\fK$ is a \define{potential kernel}{potential kernel} of $\fK$ if it satisfies (K1), (K2) above.
Below we give natural examples of potential kernels.

\begin{ex}\label{Exerighuowef}
Let $\fK$ be an arbitrary category and let $\fK_0$ consists of all monics in $\fK$.
Recall that $f$ is a \define{monic}{monic} if it satisfies the right-cancellation law:
$$f \cmp g_1 = f \cmp g_2 \implies g_1 = g_2.$$
Clearly, each identity is a monic, which proves (K1).
If $f \cmp g$ is a monic then so is $g$, showing (K2).
Note that, in general, if $f \cmp g$ is a monic then $f$ may not be a monic (take $\fK$ to be the category of sets).
This shows that the discrete norm induced by $\fK_0$ does not necessarily satisfy (\MCf).
\end{ex}

\begin{ex}
Let $\fK$ be a category of models of a fixed first order language $\El$, where the arrows are homomorphisms.
Let $\fK_0$ be the subcategory of $\fK$ whose arrows are {embeddings}.
More precisely, a homomorphism $\map f X Y$ is an \emph{embedding} if it is one-to-one and for each relation symbol $R$ in the language $\El$ it holds that
$$\seq {f(x_1), \dots, f(x_n)} \in R^Y \iff \seq {x_1, \dots, x_n} \in R^X,$$
where $R^X$, $R^Y$ are interpretations of $R$ in $X$, $Y$ respectively, and $n$ is the arity of $R$ (that is, $R^X \subs X^n$ and $R^Y \subs Y^n$).
In case where $\El$ consists of algebraic operations and constants only, every one-to-one homomorphism is an embedding.

It is clear that $\fK_0$ is a subcategory of $\fK$ with $\ob{\fK_0} = \ob{\fK}$, therefore (K1) is satisfied.
It is easy to see that $\fK_0$ satisfies (K2).
Indeed, assume $f, f \cmp g \in \fK_0$ and suppose $g \notin \fK_0$.
By Example~\ref{Exerighuowef}, we know that $g$ is one-to-one (monics in the category of sets are the one-to-one mappings).
Suppose $\map g X Y$ and $R$ is a relation symbol such that $\seq{g(x_1),\dots,g(x_n)} \in R^Y$ while $\seq {x_1, \dots, x_n} \notin R^X$.
Assume $\map f Y Z$ is an embedding.
Then $\seq {f(g(x_1), \dots, f(g(x_n))} \in R^Z$.
As $f \cmp g$ is an embedding, we conclude that $\seq {x_1, \dots, x_n} \in R^X$, a contradiction.
\end{ex}

\begin{ex}
Let $\fK$ be the category of nonempty topological spaces with reversed arrows, namely, $f \in \fK(X,Y)$ if and only if $f$ is a continuous mapping from $Y$ to $X$.
Let $\fK_0$ consist of all quotient mappings.
Recall that $\map f Y X$ is a \define{quotient mapping}{quotient mapping} if for every $U \subs X$ the set $\inv f U$ is open in $Y$ if and only if $U$ is open in $X$.
Note that $\fK_0 \subs \fK$, that is, quotient mappings are continuous.
Similar arguments as in the previous example show that $\fK_0$ is a potential kernel.
\end{ex}

\subsection{Metric spaces}\label{SubSectMetricSpacesInGeneral}

A \define{quasi-ordered category}{category!-- quasi-ordered} is a category $\fK$ in which $|\fK(a,b)| \loe 1$ for every $\fK$-objects $a,b$.
Indeed, such a category comes from a quasi-ordered set (or class, in case of large categories).
More precisely, by setting $a \loe b$ iff $\fK(a,b)\nnempty$ we obtain a reflexive and transitive binary relation on $\fK$.
A particular example of such a category is a set $X$ such that $x \loe y$ for every $x,y\in X$.
Now observe that a norm on $X$ in the sense of our definition is the same as a metric in the classical sense, except that $+\infty$ is allowed for the distance and different points can have distance zero.
Symmetry follows from Lemma~\ref{Cwergow}, because every arrow of $X$ is an isomorphism.
Note that discarding condition (\MC3) would lead to an asymmetric space, a concept studied by several authors (see~\cite{Stefan} and the references therein).
In Subsection~\ref{SubSectDigraphsAsymmetries} we explain how to deal with asymmetric spaces in our framework.

\subsection{Semigroups and groups}

Recall that a \define{semigroup}{semigroup}
is a structure of the form $\triple S \circ 1$, where $\circ$ is an associative binary operation on $S$ and $1$ is its neutral element.
A semigroup $\triple S \circ 1$ can be viewed as a category with a single object $S$, whose arrows are the elements of $S$ amd the $\circ$ is the composition. Of course, $1$ is the identity.
The conditions for a norm on $S$ are precisely (\MC1)--(\MC2). The only difference, or rather simplification, is that every two arrows are compatible.
Assume now that $S$ is a group, that is, for every $s \in S$ there is $s^{-1} \in S$ satisfying $s \circ s^{-1} = 1 = s^{-1} \cmp s$. Let $\mu$ be a norm on $S$.
Then $\mu(s^{-1}) = \mu(s)$, by Lemma~\ref{Cwergow}.
On the other hand, if $S$ is a group and $\map \mu S {[0, +\infty]}$
satisfies
\begin{enumerate}[itemsep=0pt]
	\item[(N1)] $\mu(1) = 0$.
	\item[(N2)] $\mu(f^{-1}) = \mu(f)$.
	\item[(N3)] $\mu(f \cmp g) \loe \mu(f) + \mu(g)$.
\end{enumerate}
This is indeed called a \emph{norm} on a group, leading to a left-invariant metric given by the formula $\rho(x,y) = \mu(x^{-1} \circ y)$ from which $\mu$ is trivially reconstructed.
Putting $h := f \cmp g$, we get
$$\mu(g) = \mu(f^{-1} \cmp h) \loe \mu(f^{-1}) + \mu(h) = \mu(f) + \mu(f \cmp g).$$
This shows that every function on a group satisfying (N1)--(N3) is a norm in the sense of (\MC1)--(\MC3).
Conversely, if $\mu$ is a norm on $S$ in the category-theoretic sense, then it satisfies (N1)--(N3).
Thus, our concept of a norm extends the one in group theory.

Note also that if $\pair{\fK}{\mu}$ is a normed category and $x \in \ob{\fK}$ then the group $\Aut(x)$ of all automorphisms of $x$ is a normed group, that is, it satisfies (N1)--(N3).

\subsection{A natural norm on the category of metric spaces}\label{SubSectMetsNMesur}

Below we describe one of our motivating examples of a normed category.
Namely, let $\metric$ be the category of metric spaces with Lipschitz mappings.
Let $\Lip f$ denote the Lipschitz constant of $f$, that is,
$$\Lip f = \sup_{p,q\in X,\; p\ne q}\frac{\rho_Y(f(p),f(q))}{\rho_X(p,q)},$$
where $\map f{\pair X{\rho_X}}{\pair Y{\rho_Y}}$.

Define
$$\mu(f) = \log \Bigl(\max \{ \Lip f, \Lip{\mion f} \} \Bigr),$$
where we agree that $\Lip{\mion f}=+\infty$ if $f$ is not one-to-one.
Notice that $\mu(f)=0$ if and only if $f$ is an isometric embedding.
Note also that $\mu(f) \goe 0$, as at least one of the values $\Lip f$, $\Lip{\mion f}$ is $\goe 1$.
On the other hand, if $\mu(f)=+\infty$ then $f$ may be one-to-one, although in that case its inverse is not Lipschitz.

\begin{tw}
$\pair \metric\mu$ is a normed category.
\end{tw}

\begin{pf}
It is clear that $\metric$ satisfies (\MC1). The remaining two conditions are treated below.
\begin{claim}
$\pair \metric\mu$ satisfies $(\MC2)$.
\end{claim}

\begin{pf}
Fix compatible arrows $f,g$ in $\metric$. If either $\Lip {\mion f}=+\infty$ or $\Lip {\mion g}=+\infty$ then the right-hand side of (\MC2) equals $+\infty$ and there is nothing to prove. Otherwise, we have $\Lip{f\cmp g}\loe \Lip f \cdot \Lip g$ and $\Lip{\mion{(f\cmp g)}} \loe \Lip {\mion g} \cdot \Lip {\mion f}$.
Finally, $\mu(f\cmp g)$ is either $\log \Lip{f\cmp g}$ or $\log \Lip{\mion{(f\cmp g)}}$.
\end{pf}

\begin{claim}
$\pair \metric\mu$ satisfies $(\MC3)$.
\end{claim}

\begin{pf}
Fix arrows $\map gXY$, $\map hYZ$ in $\metric$ and let $f=h\cmp g$. 
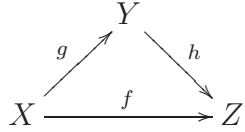
$$\xymatrix{
& Y\ar[dr]^h &\\
X\ar[ur]^g\ar[rr]^f & & Z
}$$
We shall denote the metrics on $X,Y,Z$ by the same letter $\rho$.
Let $\mu(g) = \log a$, $\mu(h) = \log b$ and $\mu(f)=\log c$.
Condition (\MC3) can be reformulated into:
\begin{equation}
a\loe bc.
\tag{1}\label{ee1}\end{equation}
In order to prove (\ref{ee1}) it suffices to check that $\rho(g(x),g(y))\loe bc \cdot \rho(x,y)$ and $\rho(g(x),g(y))\goe \mion{(bc)} \rho(x,y)$ for every $x,y\in X$.

Suppose $x,y\in X$ are such that
$\rho(g(x),g(y)) > bc \cdot \rho(x,y)$.
Then
$$c \cdot \rho(x,y) \goe \rho(f(x),f(y)) \goe \mion b \cdot \rho(g(x),g(y)) > \mion b bc \cdot \rho(x,y) = c \cdot \rho(x,y),$$
a contradiction.

Finally, suppose $x,y\in X$ are such that
$\rho(g(x),g(y)) < \mion{(bc)} \cdot \rho(x,y)$.
Then
$$\mion c \cdot \rho(x,y) \loe \rho(f(x),f(y)) \loe b \cdot \rho(g(x),g(y)) < b \cdot \mion{(bc)} \cdot \rho(x,y) = \mion c \cdot \rho(x,y),$$
which again is a contradiction.
\end{pf}
The claims above complete the proof of the theorem.
\end{pf}

Note that $\pair \metric \mu$ fails (\MCf), namely,  $|\mu(g) - \mu(h)| \loe \mu(h\cmp g)$ does not always hold.
Indeed, assume $X$ is the one-element space $\sn x$ and let $\map gXY$ be such that $Y=\dn xy$ with $\rho(x,y)=1$.
Then $\mu(g)=0$.
Now, let $Z=\dn xz$, where $\rho(x,z)=r>1$.
Let $\map h Y Z$ be the bijection satisfying $h(x)=x$.
Then $\Lip h=r>1$, so $\mu(h)>0$. On the other hand, $\mu(h\cmp g) = 0$. Hence $\mu(h)-\mu(g)\not\loe \mu(h\cmp g)$.

\subsection{Embedding-projection pairs and metric spaces}\label{SubSectionAMericaForEver}

The following example comes from America \& Rutten~\cite{America}.
The norm $\delta$ described below plays important role in the theory of reflexive domain equations developed in~\cite{America}, even though conditions (\MC1)--(\MC3) do not appear there.

Let $\epm$ denote the category whose objects are nonempty metric spaces and arrows from a metric space $X$ to a metric space $Y$ are pairs $f = \pair{e_f}{p_f}$, where $\map{e_f}{X}{Y}$ is an isometric embedding, $\map{p_f}{Y}{X}$ is a non-expansive mapping, and $p_f \cmp e_f = \id{X}$.
Given $f = \pair{e_f}{p_f}$ as above, define
$$\delta(f) = \sup_{y \in Y} \rho_Y \Bigl(e_f \cmp p_f(y), y\Bigr),$$
where $\rho_Y$ is the metric of $Y$.
Quoting from~\cite{America}: If $\map f{M_1}{M_2}$ then $\delta(f)$ can be regarded as a measure of the quality with which $M_2$ is approximated by $M_1$.

\begin{tw}
	$\delta$ is a norm on $\epm$.
\end{tw}

\begin{pf}
Clearly, $\delta(\pair{\id X}{\id X}) = 0$, therefore (\MC1) is satisfied, because the pair $\pair{\id X}{\id X}$ is the identity of $X$ in $\epm$. 

Given two mappings $f,g$ between metric spaces $X,Y$, we shall write $\rho(f,g)$ instead of $\sup_{x \in X}\rho_X(f(x), g(x))$.
Let $f = \pair{e_f}{p_f}$, $g = \pair{e_g}{p_g}$ be $\epm$-arrows for which $f \cmp g$ is defined.
Then
\begin{align*}
	\delta(f \cmp g) &= \rho(e_f \cmp e_g \cmp p_g \cmp p_f, \id{}) \loe \rho(e_f \cmp e_g \cmp p_g \cmp p_f, e_f \cmp p_f) + \rho(e_f \cmp p_f, \id{}) \\
	&\loe \rho(e_g \cmp p_g, \id{}) + \rho(e_f \cmp p_f, \id{}) = \delta(f) + \delta(g),
\end{align*}
showing that (\MC2) holds.
We have used the fact that $\rho(f \cmp g_0, f \cmp g_1) \loe \rho(g_0,g_1)$ and $\rho(g_0 \cmp h, g_1 \cmp h) \loe \rho(g_0, g_1)$. The former inequality is true because $f$ is non-expansive, the latter one holds for arbitrary mappings.

Finally, we have
\begin{align*}
	\delta(f \cmp g) + \delta(f) &= \rho(e_f \cmp e_g \cmp p_g \cmp p_f, \id{}) + \rho(e_f \cmp p_f, \id{}) \\
	&\goe \rho(e_f \cmp e_g \cmp p_g \cmp p_f, e_f \cmp p_f) \\
	&\goe \rho(e_f \cmp e_g \cmp p_g \cmp p_f \cmp e_f, e_f \cmp p_f \cmp e_f) \\
	&= \rho(e_f \cmp e_g \cmp p_g, e_f) = \rho(e_g \cmp p_g, \id{}),
\end{align*}
showing that (\MC3) holds.
Here, besides the inequalities mentioned above, we have also used the equality $\rho(e \cmp g_0, e \cmp g_1) = \rho(g_0, g_1)$, which is true whenever $e$ is an isometric embedding.
This complete the proof.
\end{pf}

\subsection{Normed linear spaces}\label{SubSectBanachNorms}

As we have already mentioned, every norm on a real or complex vector space is also a norm according to our definition, when the vector space is treated as a category.
On the other hand, the category $\nls$ of all normed linear spaces with bounded linear operators is a subcategory of the category $\metric$, because bounded linear operators are Lipschitz. Specifically, $\norm f = \Lip f$ for every linear operator $f$ between normed vector spaces.
Thus, $\nls$ is a normed category with the norm defined by
$$\mu(f) = \log \left(\norm f, \norm{f^{-1}}\right),$$
again agreeing that $\mu(f) = +\infty$ if $f$ is not 1-1 or $f^{-1}$ is unbounded.

Another possibility for a norm on $\nls$ is $\nu(f)$ defined to be the infimum of the expressions of the form $\norm{i - j \cmp f}$, where $i, j$ range through all linear isometric embeddings such that $\dom(i) = \dom(f)$, $\dom(j) = \cod(f)$, and $\cod(i) = \cod(j)$.
Note that $\nu(f)$ can be viewed as a measure of distortion that $f$ makes on its domain. We can say that the  isometric embeddings $i$, $j$ \emph{$\eps$-correct} $f$ if $\norm{i - j \cmp f} \loe \eps$. Thus, $\nu(f)$ measures how difficult it is to correct $f$.
The fact that $\nu$ is a norm (proved below) relies on the \define{amalgamation property}{amalgamation property} of $\nls$, namely, for every $\nls$-arrows $\map f Z X$, $\map g Z Y$ there exist $\nls$-arrows $\map{f'}{X}{W}$, $\map{g'}{Y}{W}$ such that $f' \cmp f = g' \cmp g$. Furthermore, if $f$, $g$ are isometric embeddings, then $f'$, $g'$ could be chosen to be isometric embeddings, too.
The amalgamation property (belonging to the folklore of Banach space theory) is an important tool for constructing special Banach spaces, see~\cite{ACCGM, GarKubExtracta}. 

\begin{tw}
	$\pair{\nls}{\nu}$ is a normed category.
\end{tw}

\begin{pf}
	Clearly, $\nu(\id X) = 0$, therefore (\MC1) holds.
	Fix $\map f Y Z$, $\map g X Y$ in $\nls$, so that $f \cmp g$ is defined.
	
	Suppose $\nu(f \cmp g) > \nu(f) + \nu(g)$ and choose linear isometric embeddings $i_f$, $j_f$, $i_g$, $j_g$ such that
	$$\nu(f \cmp g) > \norm{i_f - j_f \cmp f} + \norm{i_g - j_g \cmp g}.$$
	Let $k$, $\ell$ be an amalgamation of $j_g$, $i_f$, that is, $k$, $\ell$ are linear isometric embeddings such that $k \cmp j_g = \ell \cmp i_f$.
	Now
	\begin{align*}
		\nu(f \cmp g) &\loe \norm{k \cmp i_g - \ell \cmp j_f \cmp f \cmp g} \\
		&\loe \norm{k \cmp i_g - k \cmp j_g \cmp g} + \norm{k \cmp j_g \cmp g - \ell \cmp j_f \cmp f \cmp g} \\
		&\loe \norm{i_g - j_g \cmp g} + \norm{\ell \cmp i_f \cmp g - \ell \cmp j_f \cmp f \cmp g} \\
		&\loe \norm{i_g - j_g \cmp g} + \norm{i_f - j_f \cmp f} < \nu(f \cmp g),
	\end{align*}
	a contradiction.
	This shows (\MC2).
	
	Now suppose that $\nu(g) > \nu(f \cmp g) + \nu(f)$ and choose $i_f, j_f, i_{fg}, j_{fg}$ witnessing this fact. Let $k, \ell$ be linear isometric embeddings providing the amalgamation of $j_f, j_{fg}$, that is, $k \cmp j_f = \ell \cmp j_{fg}$. We have
	\begin{align*}
		\nu(g) &\loe \norm{\ell \cmp i_{fg} - k \cmp i_f \cmp g} \\
		&\loe \norm{\ell \cmp i_{fg} - \ell \cmp j_{fg} \cmp f \cmp g} + \norm{\ell \cmp j_{fg} \cmp f \cmp g - k \cmp i_f \cmp g} \\
		&=\norm{\ell \cmp i_{fg} - \ell \cmp j_{fg} \cmp f \cmp g} + \norm{k \cmp j_{f} \cmp f \cmp g - k \cmp i_f \cmp g} \\
		&\loe \norm{i_{fg} - j_{fg} \cmp f \cmp g} + \norm{i_f - j_f \cmp f} < \nu(g),
	\end{align*}
	a contradiction. This shows (\MC3) and completes the proof.
\end{pf}

The norm $\nu$ appears implicitly in~\cite{KubSol}, playing a crucial role there. It turns out that $\nu$ has a nice correlation with the usual norm:

\begin{tw}\label{THMsbkfiaofj}
	Let $\map f X Y$ be a linear operator between normed spaces, $0 < \eps < 1$.
	The following conditions are equivalent.
	\begin{enumerate}[itemsep=0pt]
		\item[{\rm(a)}] There exists a norm on $X \oplus Y$ such that the canonical embeddings $\map{e_X}{X}{X \oplus Y}$, $\map{e_Y}{Y}{X \oplus Y}$ are isometric and $\norm{e_X - e_Y \cmp f} \loe \eps$.
		\item[{\rm(b)}] $(1 - \eps) \norm x \loe \norm{f(x)} \loe (1 + \eps) \norm x$ holds for every $x \in X$.
	\end{enumerate}	
\end{tw}

\begin{pf}
	(a)$\implies$(b)
	Fix $x \in X$ with $\norm x = 1$.
	We have
	$$\norm{f(x)} = \norm{e_Y f(x)} \loe \norm{e_Y f(x) - e_X(x)} + \norm{e_X(x)} \loe \eps + 1$$
	and
	$$\norm{f(x)} = \norm{e_Y f(x)} \goe \norm{e_X(x)} - \norm{e_Yf(x) - e_X(x)} \goe 1 - \eps.$$
	(b)$\implies$(a)
	A suitable norm on $X \oplus Y$ can be defined by the formula
	$$\norm{\pair xy} = \inf \bigsetof{\norm{x-w}_X + \norm{y+f(w)}_Y + \eps \norm{w}_X}{w \in X},$$
	where $\anorm_X$ and $\anorm_Y$ are the norms of $X$ and $Y$, respectively. Checking that (a) holds is easy, yet it requires some computations. We refer to~\cite[Proof of Lemma 3.1]{CGK} for details, where it is done for the case of $p$-norms, where $0<p \loe 1$.
\end{pf}

The norm $\nu$ can also be defined on the category of metric spaces, where Theorem~\ref{THMsbkfiaofj} would still hold, as there is a functor embedding isometrically each metric space into a normed linear space.
The norm $\nu$ is actually a special case of a much more abstract construction in categories enriched over metric spaces, see~\cite{Kub41}.

\subsection{Directed graphs, free categories, asymmetric spaces}\label{SubSectDigraphsAsymmetries}

Recall that a \define{weight}{weight} on a graph (directed or undirected) is a function assigning to each arrow/edge a nonnegative real value.
A weighted undirected graph naturally leads to a metric space.
Namely, the distance between two vertices $a$, $b$ is declared to be the infimum of the weights of all possible paths between $a$ and $b$, assuming that it is $+\infty$ if $a$ and $b$ are in different components of the graph.
If the graph is infinite, this might be just a pseudo-metric, as the value zero may be obtained for two different vertices.
Recall that the weight of a path is the sum of the weights of all its arrows/edges.

It is clear how to define a quasi-metric (that is, a function fulfilling all the axioms of a metric except the symmetry) on an arbitrary directed graph.

To be more precise, recall that a \define{digraph}{digraph} (or a \define{directed graph}{digraph}) is a quadruple $\seq{ V, A, \dom, \cod }$,
where $\map \dom A V$ and $\map \cod A V$ are functions, often called the \emph{source} and the \emph{sink}, respectively.
There are no rules (axioms) here.
The elements of $A$ are called \emph{arrows}, while the elements of $V$ are called \emph{vertices} or \emph{objects}.
In this setting, a \define{category}{category} is nothing but a digraph endowed with a composition operation satisfying the usual axioms, adding the requirement that each vertex has a distinguished loop, called \emph{identity}.
Given a digraph $X = \seq{V, A, \dom, \cod}$, the \define{free category}{free category} $\Fr(X)$ over $X$ is defined to be the collection of all triples $\triple x p y$, where $x, y \in V$ and $p$ is a path from $x$ to $y$, understood in the usual way.
The composition is defined in the obvious way.
The empty set is considered to be a path from an arbitrary vertex to itself, therefore $\triple x \emptyset x$ is the identity of $x \in V$.
In particular, $V$ is the set of objects of $\Fr(X)$.

If $X$ is a weighted digraph with a weight function $\map w A {[0,+\infty)}$, the free category $\Fr(X)$ has a natural norm given by
$$\mu(x,p,y) = \sum_{i=1}^n w(f_i),$$
where $p = \seq{f_1,f_2, \dots, f_n}$.
Note that $\mu$ satisfies a much stronger condition than (\MC2) and (\MC3), namely:
$$\mu(f \cmp g) = \mu(f) + \mu(g)$$
whenever $f \cmp g$ is defined.
In particular, $\mu(f) \loe \mu(f \cmp g)$ and $\mu(g) \loe \mu(f \cmp g)$, which is an evidence even for (\MCf).

Finally, we would like to represent asymmetric spaces as normed categories.
For this aim, let $\pair X \rho$ be an asymmetric space and let us convert it to a digraph in which for each $x,y \in X$ there is a unique arrow $\overrightarrow{x y}$.
The paths in $X$ are of the form
$$p = \seq{\overrightarrow{x_1 x_2}, \overrightarrow{x_2 x_3}, \dots, \overrightarrow{x_{n-1}x_n}}$$
and we have $\mu(p) = \rho(x_1,x_2)+\rho(x_2,x_3)+ \dots +\rho(x_{n-1},x_n)$.
By this way we obtain a norm $\mu$ on $\Fr(X)$. Note that the original quasi-metric $\rho$ is reconstructible as $\mu$ restricted to paths of length $1$.
This approach suggests that for some purposes it might be meaningful to consider directed graphs instead of quasi-metric spaces.

\subsection{Quasi-metrics from normed categories}

Actually, there is also a way of obtaining a quasi-metric from a normed category $\pair \fK \mu$.
Namely, let $X := \ob \fK$ be the class of all objects of $\fK$.
Given $x, y \in X$, define
$$\rho(x,y) = \inf_{f \in \fK(x,y)} \mu(f),$$
where $\fK(x,y)$ denotes (as usual) the set of all $\fK$-arrows from $x$ to $y$.
We agree that $\rho(x,y) = +\infty$ if $\fK(x,y) = \emptyset$.
Conditions (\MC1), (\MC2) guarantee that $\rho$ satisfies
\begin{enumerate}
\item[(Q1)] $\rho(x,x) = 0$,
\item[(Q2)] $\rho(x,y) \loe \rho(x,z) + \rho(z,y)$,
\end{enumerate}
which is the standard definition of a quasi-metric, allowing $+\infty$ as a possible value.
This might serve as a motivation for further study in the theory of asymmetric spaces, although this line of research will not be explored here.

\separator

We believe that the examples described above provide enough justification for considering our set of axioms (\MC1)--(\MC3), while at the same time discarding (\MCf).

\section{Cauchy completion}\label{SectCauchySeqsCmps}

Let $\pair \fK \mu$ be a fixed normed category.
A \define{sequence}{sequence} in $\fK$ is, by definition, a covariant functor from $\omega$ (the set of natural numbers) into $\fK$.
We shall denote sequences by $\wek x, \wek y$, etc.
Given a sequence $\map{\wek x}\omega \fK$, we denote by $x_n$ the $n$th element of $\wek x$ and by $x_n^m$ the bonding arrow from $x_n$ to $x_m$.

\begin{df}
	A sequence $\map {\wek x}\omega \fK$ is \define{convergent}{sequence!-- convergent} if there exist $v\in \ob \fK$ and a family of arrows $\sett{\map{g_n}{x_n}v}\ntr$ such that
	\begin{enumerate}[itemsep=0pt]
		\item[(C1)] $\pair v{\ciag{g}}$ is the co-limit of $\wek x$, i.e. $\map{g_n}{x_n} v$ are such that $g_m\cmp x_n^m=g_n$ for $n<m$ and for every pair $\pair {v'}{\sett{g'_n}{\ntr}}$ with the same property there exists a unique arrow $\map h v {v'}$ satisfying $h\cmp g_n = g_n'$ for $\ntr$.
		\item[(C2)] $\lim_{n\to\infty} \mu(g_n) = 0$.
	\end{enumerate}
	We shall write $\pair v{\ciag g}=\lim\wek x$ or just $v = \lim\wek x$ when there is no confusion what the arrows $g_n$ are taken into account.
	The object $v$ (more formally, the pair $\pair v{\ciag g}$) is the \define{limit}{sequence!-- limit} of $\wek x$.	
\end{df}

Unfortunately, there is no way to prove that the limit is unique up to a $0$-isomorphism.
Here is an example:

\begin{ex}\label{Exeihgouw}
Let $\fK$ be a poset category whose class of objects is $\omega \cup \dn{\ell_0}{\ell_1}$, where $\omega$ has the usual ordering and $n \loe \ell_i$ for every $\ntr$, $i = 0,1$.
We make $\ell_0$ and $\ell_1$ isomorphic, setting $\ell_0 \loe \ell_1$ and $\ell_1 \loe \ell_0$.
By this way, $\loe$ is a quasi-ordering on the $\fK$-objects.
We now define $\mu(x,y) = 0$ whenever $x \loe y$ and $x \notin \dn{\ell_0}{\ell_1}$.
Finally, we set $\mu(\ell_0, \ell_1) = +\infty = \mu(\ell_1, \ell_0)$.
By this way, $\pair \fK \mu$ becomes a normed category.
Let $\wek x$ be the sequence such that $x_n = n$ for $\ntr$.
Then $\lim \wek x = \ell_i$ for $i = 0,1$.
Of course, $\ell_0$ and $\ell_1$ are isomorphic, however they are not $0$-isomorphic.
\end{ex}

\begin{df}
	A \define{Cauchy sequence}{sequence!-- Cauchy} is a sequence $\wek x$ satisfying the following condition:
	\begin{enumerate}
		\item[(C)] For every $\eps > 0$ there exists $n_0$ such that $\mu(x_m^n) < \eps$ whenever $n_0 < m \loe n$.
	\end{enumerate}
\end{df}
Particular case of a Cauchy sequence is any sequence $\wek x$ such that $\mu(x_k^\ell)=0$ for every $k < \ell < \nat$. Such a sequence will be called \define{stable}{sequence!-- stable}.
As in the case of metric spaces, one can use convergent series to test whether a sequence is Cauchy.
Namely, we have the following simple fact, which is a consequence of (\MC2):

\begin{prop}\label{PropSeriesConvrgnt}
Let $\pair \fK \mu$ be a normed category and let $\wek x$ be a sequence in $\fK$.
If $\sum_{n=0}^\infty \mu(x_n^{n+1}) < +\infty$ then $\wek x$ is a Cauchy sequence.
\end{prop}

One needs to be careful here, namely, if a subsequence is Cauchy then the sequence may not be Cauchy (of course, this pathology does not occur in metric spaces).
Here is an example:

\begin{ex}
Let $\Nat$ be the set of natural numbers, treated as a poset category.
Define $\mu(m,n) = 0$ if $m$ is even or $m=n$, and $\mu(m,n) = +\infty$, otherwise (here $m < n$).
Clearly, $\pair \Nat \mu$ is a normed category.
The natural sequence $\wek x$, $x_n = n$, is not Cauchy, although its subsequence obtained by restricting the domain to all even numbers is Cauchy.
\end{ex}

Axiom (\MC3) is needed to yield the following natural fact.

\begin{lm}
Every convergent sequence is Cauchy.
\end{lm}

\begin{pf}
Let $\wek x$ be a sequence convergent to $\pair g{\sett{g_n}{\ntr}}$. Fix $\eps>0$ and find $n_0$ such that $\mu(g_n)<\eps/2$ whenever $n\goe n_0$. Given $n<m$, using (\MC3) we get
$$\mu(x^m_n) \loe \mu(g_m\cmp x^m_n) + \mu(g_m) = \mu(g_n) + \mu(g_m)<\eps.$$
Thus $\mu(x^m_n)<\eps$ whenever $m>n\goe n_0$.
\end{pf}

We shall denote by $\ciagi(\fK, \mu)$ the class of all Cauchy sequences in $\pair \fK \mu$.
Given $\wek x, \wek y \in \ciagi(\fK)$, a \define{transformation}{sequence!-- transformation} from $\wek x$ to $\wek y$ is a natural transformation from $\wek x$ to $\wek y \cmp \phi$, where $\map \phi \omega \omega$ is increasing (i.e. a self-functor on $\omega$).
More precisely, a transformation is a family of arrows $\wek f  = \ciag f$ such that $\map {f_n}{x_n}{y_{\phi(n)}}$ for some increasing map $\map \phi \omega \omega$, satisfying
$$y_{\phi(m)}^{\phi(n)} \cmp f_m = f_n \cmp x_m^n$$
for every $m < n < \omega$.
We shall write $\map {\wek f}{\wek x}{\wek y}$.

It is clear how to identify a normed category $\pair \fK \mu$ with a subcategory of $\ciagi(\fK, \mu)$; namely, every $\fK$-object $x$ corresponds to the constant sequence
$$\xymatrix{x \ar[r] & x \ar[r] & x \ar[r] & \dots}$$
where all the arrows are identities.
Further on, we shall make this identification.
Note that every transformation between constant sequences is induced by a uniquely determined $\fK$-arrow.
In order to make the theory reasonable, it remains to show that $\ciagi (\fK, \mu)$ is again a normed category and that all Cauchy sequences converge in $\ciagi (\fK, \mu)$.
This is done below.

\begin{lm}\label{Lierhif}
Let $\map {\wek f}{\wek x}{\wek y}$ be a 
transformation of Cauchy sequences. Then the limit
$$\lim_{n \to \infty}\mu(f_n)$$
exists.
\end{lm}

\begin{pf}
Assume $\map {f_n}{x_n}{y_{\phi(n)}}$, where $\map \phi \omega \omega$ is increasing.
Let $$\ell = \liminf_{n\to \infty}\mu(f_n)$$ and let $\sett{k_n}{\ntr}$ be such that $\lim_{n\to \infty}\mu(f_{k_n}) = \ell$.
We may assume that $n \loe k_n$ for every $\ntr$.
Note that
$$f_m \cmp x_n^m = y_{\phi(n)}^{\phi(m)} \cmp f_n$$ whenever $n \loe m$.
Using this equation with $m=k_n$ and applying (\MC2), (\MC3), we get
\begin{equation}
\mu(f_n) \loe \mu(y_{\phi(n)}^{\phi(k_n)} \cmp f_n) + \mu(y_{\phi(n)}^{\phi(k_n)}) \loe \mu(x_n^{k_n}) + \mu(f_{k_n}) + \mu(y_{\phi(n)}^{\phi(k_n)}).
\tag{*}\label{eqsegorosgo}
\end{equation}
Note that
$$\lim_{n \to \infty} \mu(x_n^{k_n}) = 0 \oraz \lim_{n \to \infty} \mu(y_{\phi(n)}^{\phi(k_n)}) = 0,$$
because $\wek x$ and $\wek y$ are Cauchy sequences.
Using this together with (\ref{eqsegorosgo}), we obtain
$$\limsup_{n \to \infty}\mu(f_n) \loe \ell.$$
It follows that $\ell = \lim_{n \to \infty}\mu(f_n)$.
\end{pf}

By the lemma above, we conclude that $\ciagi(\fK,\mu)$ is again a normed category with the norm defined by
$$\ovr\mu(\wek f) := \lim_{n \to \infty} \mu(f_n).$$
It is easy to check that conditions (\MC1)--(\MC3) are satisfied.
Note also that $\mu(\wek f) = \mu(f)$ whenever $\wek f$ is a transformation between sequences of identities determined by a $\fK$-arrow $f$.

A normed category $\pair \fK \mu$ is said to be \define{Cauchy complete}{category!-- Cauchy complete} if every Cauchy sequence is convergent in $\pair \fK \mu$.
It is natural to expect that $\ciagi(\fK,\mu)$ is Cauchy complete. This is indeed the case, as we show below.

\begin{tw}
Let $\pair \fK \mu$ be a normed category.
Then $\ciagi(\fK,\mu)$ is Cauchy complete.
\end{tw}

\begin{pf}
This is a standard diagonalization argument.
Namely, let $\wek {\wek x}$ be a Cauchy sequence in $\ciagi(\fK,\mu)$.
We make the agreement that $x_{n,m}$ will denote the $m$th element of the sequence $\wek{x}_n$.
Refining the sequences $\wek{x}_n$ inductively, we may assume that the transformations from $\wek{x}_n$ to $\wek{x}_{n+1}$ are natural.
That is, the sequence $\wek {\wek x}$ can be presented as the following infinite matrix:
$$\xymatrix{
\vdots & \vdots & \vdots & \\
x_{0,2} \ar[r] \ar[u] & x_{1,2} \ar[r] \ar[u] & x_{2,2} \ar[r] \ar[u] & \cdots \\
x_{0,1} \ar[r] \ar[u] & x_{1,1} \ar[r] \ar[u] & x_{2,1} \ar[r] \ar[u] & \cdots \\
x_{0,0} \ar[r] \ar[u] & x_{1,0} \ar[r] \ar[u] & x_{2,0} \ar[r] \ar[u] & \cdots
}$$
Notice that the diagonal sequence
$$\xymatrix{
x_{0,0} \ar[r] & x_{1,1} \ar[r] & x_{2,2} \ar[r] & \cdots
}$$
is the co-limit of $\wek {\wek x}$.
The only obstacle is that this sequence may not be Cauchy and (C2) may not be satisfied.
However, this can be achieved by further refining the sequences $\wek{x}_0, \wek{x}_1, \wek{x}_2, \dots$.
\end{pf}

A natural way to compare normed categories is by using `norm-preserving' functors.
Let $\pair \fK \mu$, $\pair \fL \nu$ be fixed normed categories.
A functor $\map F \fK \fL$ is \define{non-expansive}{functor!-- non-expansive} if
$$\nu(F(f)) \loe \mu(f)$$
for every $\fK$-arrow $f$.
This is an obvious generalization of non-expansive maps between metric spaces.
A functor is \define{sequentially continuous}{functor!-- sequentially continuous} if it preserves limits of convergent sequences.

It turns out that $\ciagi(\fK,\mu)$ has the following universal property:

\begin{tw}
Let $\pair \fK \mu$ be a normed category.
Given a Cauchy complete normed category $\pair \fL \nu$, for every sequentially continuous non-expansive functor $\map F \fK \fL$ there exists a unique sequentially continuous functor $\map {\ovr F}{\ciagi(\fK,\mu)} \fL$ extending $F$.
\end{tw}

\begin{pf}
Put $\ovr F(\wek x) = \lim  (F \cmp \wek{x})$.
Clearly, this defines a functor from $\ciagi(\fK,\mu)$ to $\fL$.
Uniqueness of $\ovr F$ is obvious.
It remains to check that $\ovr F$ is indeed sequentially continuous.
However, this follows from the following property of co-limits:
$$\lim \wek{\wek x} = \lim_{n\to\infty} \lim_{m\to\infty} x_{n,m}.$$
\end{pf}

We cannot prove that $\ovr F$ is non-expansive, see Example~\ref{Exeihgouw} above.
Actually, uniqueness of $\ovr F$ is up to an isomorphism only, not necessarily a 0-isomorphism.

In view of the statement above, $\ciagi(\fK,\mu)$ will be called the \define{Cauchy completion}{Cauchy completion} of $\pair \fK \mu$.

\separator

After all, we do not care much about the uniqueness of limits of convergent sequences.
In natural examples the limits are indeed unique.
On the other hand, we never rely on limit uniqueness; the only problem is to make the theory more natural and aesthetic...
Anyway, it makes sense to call a normed category $\pair \fK \mu$ \define{Hausdorff}{normed category!-- Hausdorff} if every converging sequence in $\fK$ has a unique limit, up to $0$-isomorphisms.

Below are examples of a Cauchy complete category.

\begin{ex}
	Let $\metric$ be the category of metric spaces endowed with the usual norm (see Subsection~\ref{SubSectMetsNMesur}). Let $\vec x$ be a Cauchy sequence in $\metric$. Removing its initial part, we may assume that the bonding mappings $\map{x_n^m}{X_n}{X_m}$ are one-to-one, therefore without loss of generality we may even assume that they are inclusions. Thus $X_n \subs X_m$ for $n \loe m$ and the only problem is that the metric $\rho_n$ on $X_n$ may not coincide with $\rho_m$.
	Let $X_\infty = \bigcup_{\ntr} X_n$.
	The fact that $\vec x$ is Cauchy implies that the limit
	$$\rho_\infty(s,t) := \lim_{n \to \infty} \rho_n(s,t)$$
	exists for every $s,t \in X_\infty$. Furthermore, $\rho_\infty$ is obviously a metric. It is rather obvious that $X_\infty$ is the limit of $\vec x$.
	
	Another, perhaps more useful, example is the category $\ovr \metric$ of all complete metric spaces, endowed with the same norm.
	It is Cauchy complete too, as one can take $X_\infty$ to be the completion of $\pair{\bigcup_{\ntr}}{\rho_\infty}$.
	This construction of a limit is particularly useful when a metric (or a norm, in the case of Banach spaces) needs to be defined by certain approximations.
\end{ex}

\begin{ex}
	Let $\pair \epm \delta$ be the category of embedding-projection pairs considered in~\cite{America} and discussed in Subsection~\ref{SubSectionAMericaForEver}.
	It has been proved in~\cite{America} that $\pair \epm \delta$ is Cauchy complete.
	It is worth pointing out that arbitrary sequences in $\epm$ do not necessarily have co-limits, the norm $\delta$ is essential to get convergence.
\end{ex}

\section{Banach's fixed-point theorem}\label{SectBanachFixedpoint}

There is no need to recall the famous Banach's Contraction Principle.
On the other hand, it is hopeless to try to create a comprehensive list of its generalizations, scattered through the literature.
Still, once we have made a proper generalization of the concept of a complete metric space, it has to be tested on the validity of Banach's fixed-point theorem.
It turns out that such a statement is not only true in our setting, but its proof becomes merely trivial!
Triviality is actually the consequence of our setting: we work with Cauchy sequences as the main objects of investigations.

Given two normed categories, it is clear how to define the notion of a \define{contraction}{contraction}: it is a covariant functor $F$ satisfying $\mu(F(f)) \loe r \mu(f)$ for every arrow $f$, where $r < 1$ is fixed and independent of $f$.
A \define{fixed-point}{fixed-point} of $F$ is an object $X$ for which there exists a $0$-isomorphism $\map{h}{X}{F(x)}$.
The classical Banach's Contraction Principle can be stated as follows:

\begin{tw}\label{ThmBCPrevisited}
Let $\pair \fK \mu$ be a normed category and let $\map F \fK \fK$ be a contraction.
Suppose $\map f a {F(a)}$ is a $\fK$-arrow such that $\mu(f) < +\infty$.
Then the sequence
\begin{equation}
\xymatrix{
a \ar[r]^f & F(a) \ar[r]^{F(f)} & F^2(a) \ar[r]^{F^2(f)} & F^3(a) \ar[r]^{F^3(f)} & \cdots
}
\tag{$\star$}\label{EqSeqFixdPt}
\end{equation}
is a fixed-point of the extension $\ovr F$ defined on $\ciagi(\fK,\mu)$.
\end{tw}

The proof is rather obvious: Apply the functor $F$ to the sequence (\ref{EqSeqFixdPt}) and draw a suitable infinite diagram in order to see that $\lim_{n \to \infty} F^n(f)$ is a $0$-isomorphism. The fact that $F$ is a contraction implies that the sequence (\ref{EqSeqFixdPt}) is Cauchy.

\begin{wn}
	Assume $\pair \fK \mu$ is a Cauchy complete normed category, $\map F \fK \fK$ is a sequentially continuous contraction. If there exists a $\fK$-arrow of the form $\map f a {F(a)}$ with $\mu(f)$ finite, then $F$ has a fixed-point, i.e., there exist $s \in \ob{\fK}$ and a $0$-isomorphism $\map h s {F(s)}$.
\end{wn}

Theorem~\ref{ThmBCPrevisited} indeed contains the essence of the classical Banach's fixed-point theorem. This is evidenced by the case where $\fK$ is a metric space, viewed as a normed category, when there is always a unique arrow
$$a \to F(a)$$
for every $a \in \ob \fK$.
Note that Theorem~\ref{ThmBCPrevisited} says in particular that the sequence (\ref{EqSeqFixdPt}) is Cauchy, which follows from Proposition~\ref{PropSeriesConvrgnt}, knowing from basic calculus that the series $\sum_{n=1}^{+\infty}r^n$ is convergent whenever $r<1$.

Special cases of Theorem~\ref{ThmBCPrevisited} already appeared in the category-theoretic literature, typically in theoretical computer science, where one of the lines of research is solving equations of the form $X \iso F(X)$, where $F$ is some functor.
It seems that America \& Rutten~\cite{America} is the first work in this topic containing explicitly formulated fixed-point theorem.
Theorem~\ref{ThmBCPrevisited} says nothing about uniqueness of fixed-points, as in general it is impossible to get it:

\begin{ex}
	Let $\pair \fK \mu$ be a Cauchy complete category with a sequentially continuous contraction $\map F \fK \fK$.
	Let $\fK \oplus \fK$ be the disjoint sum of two copies of $\fK$ and let $\mu \oplus \mu$ and $F \oplus F$ be defined in the obvious way. It is clear that $F \oplus F$ is a contraction and each fixed-point of $F$ appears now in two non-isomorphic copies that are fixed-points of $F \oplus F$.
\end{ex}

Thus, in order to get uniqueness of fixed-points (which is trivial in metric spaces) one needs to look for some reasonable extra conditions. This issue is discussed in~\cite{America} and in several subsequent works, including~\cite{Birkedal} where categories enriched over ultrametric spaces were considered.

\section{Final remarks}\label{SectFinals}

We believe that the concept of a norm satisfying (\MC1)--(\MC3) is applicable and worth further studying.
We have shown that several natural categories have canonical norms.
Some of them were essentially used, see e.g.~\cite{America, Birkedal, Grandis, Kent, Macaj, Zlatos}, where the authors often discard (\MC3).

Our personal motivation for introducing norm on categories comes from our preprint \cite{Kub41}, where one deals with a metric-enriched category $\fK$ with a fixed subcategory $\fK_0$. The norm of $f \in \fK$ is defined to be the infimum of all distances $\rho(i, j \cmp f)$, where $i,j \in \fK_0$ are such that $\dom(i) = \dom(f)$, $\dom(j) = \cod(f)$, and $\cod(i) = \cod(j)$. It turns out that $\mu$ satisfies (\MC1)--(\MC3) provided that $\ob{\fK_0} = \ob{\fK}$ and $\fK_0$ has the \emph{almost amalgamation property}, that is, for every $i,j \in \fK_0$ with the same domain, for every $\eps>0$ there exist $i', j' \in \fK_0$ such that $\dom(i') = \cod(i)$, $\dom(j') = \cod(j)$, $\cod(i') = \cod(j')$, and $\rho(i' \cmp i, j' \cmp j) < \eps$.
This norm is needed for defining the concept of a \emph{generic sequence} which is crucial in studying universal almost homogeneous objects, also called \emph{continuous \fra\ limits}.

Concerning further research,
we believe that several well-established concepts (e.g. sequential compactness) will find their place in the new, more general, context.
In fact, more detailed study of abstract convergence involving Cauchy sequences (that clearly require norms) had already been indicated in the last section of America \& Rutten~\cite{America}.
The main point of further research studying norms is of course finding non-trivial results and applications, at the same time avoiding trivialities.

\printindex

\end{document}